\documentclass{amsart}
\usepackage{amsmath}
\usepackage{amsthm}
\usepackage{amsfonts}
\usepackage{amsopn}
\usepackage{amssymb}
\usepackage{rlepsf}
\usepackage[matrix,arrow,curve,cmtip]{xy}

\newcommand{\R}{\ensuremath{\mathbb{R}}}

\newcommand{\N}{\ensuremath{\mathbb{N}}}
\newcommand{\Z}{\ensuremath{\mathbb{Z}}}

\DeclareMathOperator{\image}{image}

\newtheorem{thm}{Theorem}
\newtheorem{prop}[thm]{Proposition}
\newtheorem{lemma}[thm]{Lemma}
\begin{document}
\title{Averaged Dehn Functions for Nilpotent Groups}
\author{Robert Young }
\address{Institut des Hautes \'Etudes Scientifiques, Le Bois Marie, 35 route de Chartres F-91440 Bures-sur-Yvette, France}

\thanks{The author gratefully acknowledges support from a GAANN
  fellowship during part of the writing of this paper.}
\bibliographystyle{plain}
\begin{abstract}
Gromov proposed an averaged version of the Dehn function and claimed
that in many cases it should be subasymptotic to the Dehn function.
Using results on random walks in nilpotent groups, we confirm this
claim for most nilpotent groups.  In particular, if a nilpotent group
satisfies the isoperimetric inequality $\delta(l)<Cl^\alpha$ for
$\alpha>2$ then it satisfies the averaged isoperimetric inequality
$\delta^{\text{avg}}(l)<C'l^{\alpha/2}$.  In the case of non-abelian free
nilpotent groups, the bounds we give are asymptotically sharp.
\end{abstract}
\maketitle
\section{Introduction}

Determining the asymptotic behavior of the Dehn function for a group
is a much-studied problem in group theory; see \cite{Bri} for an
introduction to the subject.  In \cite[5.$\text{A}'_6(\text{c})$, p. 90]{Gro}, Gromov proposed a
variation on this problem: instead of asking for the largest area
required to fill a closed curve of a given length, he asked what the
average area is, taken over all closed curves.  Gromov claimed that in
many cases, this averaged Dehn function should be asymptotically
smaller than the Dehn function, which was confirmed in the case of
finite rank free abelian groups by \cite{KukRom}, whose bound was
improved by Bogopolski and Ventura \cite{BogVen}.  Still, little is
known about the averaged Dehn function; some open questions include
whether it is invariant under changes of generators or
quasi-isometries.

In this paper, we prove upper and lower bounds for the averaged
Dehn function of a nilpotent group which show, in particular, that if
$G$ is a nilpotent group which isn't finite or virtually $\Z$, then
its averaged Dehn function is subasymptotic to its Dehn function.
This implies that random walks in nilpotent groups bound much less
area on average than the worst-case curves.  In addition, our bounds
are sharp in many cases.

A related problem appears in statistical mechanics, where the behavior
of a charged particle moving randomly in a magnetic field depends on
the signed area of its path.  The distributions of the signed area and
the area of a random loop have thus been considered by mathematical
physicists, see for instance, \cite{Colomo} and \cite{BeCaBaCl}.  In
the former, Colomo shows that our bounds break down for infinitely
generated groups.  He considers the area of a random loop in $\Z^d$
and shows that if $n$ is fixed and $d\to\infty$,
$\delta^{\text{avg}}(n)\to n(n-1)/6$, in contrast to the finite
dimensional case, where $\delta^{\text{avg}}(n)$ grows strictly
subquadratically.

The averaged Dehn function can also be interpreted as reflecting
properties of the average-case complexity of the word problem for a
group just as the Dehn function reflects its complexity.  Kapovich,
Myasnikov, Schupp, and Shpilrain studied this average-case complexity
by showing that in many groups, it is easy to show that most elements
are not the identity \cite{KaMyScSh}; the averaged Dehn
function represents the complexity of verifying that an element
represents the identity.  One interesting question, then, might be to
find groups where the averaged Dehn function differs substantially
from the Dehn function; the results in this paper show that for many
nilpotent groups, $\delta_{\text{avg}}(n)$ is approximately
$\delta_{\text{avg}}(\sqrt{n})$.

In section \ref{defs}, we define Dehn functions and averaged Dehn
functions and state the upper bound.  In section \ref{centiso}, we
define the centralized isoperimetric function of a group, give a
method to calculate it for a nilpotent group, and state a lower bound
using this function.  In sections \ref{lowbndsec} and \ref{upbndsec},
we prove these bounds.

The author would like to thank Shmuel Weinberger and Benson Farb for
their helpful suggestions and to thank the referee for their detailed
suggestions on improving this paper.  The work in this paper was done
as part of the author's doctoral thesis at the University of Chicago.

\section{Definitions}\label{defs}
We will be using big O notation for asymptotic bounds throughout this
paper.  Recall that $O$ represents an asymptotic upper bound, $\Omega$
represents an asymptotic lower bound, and $\Theta$ represents an
asymptotically tight bound.  Specifically,
\begin{align*}
f(x)=O(g(x)) & \text{iff $\exists M, x_0$ s.t. $f(x)\le Mg(x)$ for $x>x_0$}\\
f(x)=\Omega(g(x)) & \text{iff $\exists m, x_0$ s.t. $mg(x)\le f(x)$ for $x>x_0$}\\
f(x)=\Theta(g(x)) & \text{iff $\exists m,M, x_0$ s.t. $mg(x)\le f(x)\le Mg(x)$ for $x>x_0$}.
\end{align*}

One can define the Dehn function in a variety of contexts; here, we
define it in
terms of a presentation of a group.  We first define the
filling area of a word.

Let $G$ be a group with identity $e$, given by a presentation
$$G=\{e_1,\dots,e_d|r_1,\dots,r_s\}.$$  Let $R=\langle r_1,\dots,r_s\rangle$ be the
normal closure of the relators, and $F$ be the free group on $d$ generators, so that $G=F/R$.  If $w$ is a word
in the $e_i^{\pm 1}$ that is the identity in $G$, then $w$ lies in $R$ when considered as an element of $F$.  We can thus write
$$w=\prod_{i=1}^k g_i^{-1} r_{a_i}^{\pm1} g_i,$$
where the equality is taken in $F$.  Define $\delta_G(w)$ as the minimal $k$ for which we can write such a
decomposition.  We will drop the $G$ if context makes it clear which group is meant.

$\delta(w)$ counts the number of applications of relators required to
reduce $w$ to the trivial word.  We can view this as an area by
considering the $2$-complex obtained by taking the Cayley graph of $G$
and adding a face for each conjugate of a relator.  Then $w$
represents a curve in this complex, and $\delta(w)$ represents the
minimal number of faces in a disc with boundary $w$.  Alternately, we
can view $\delta(w)$ as a metric on $R$; it is the word metric on $R$
corresponding to the (usually infinite) generating set consisting of all conjugates of
the $r_i$.

We define the Dehn function of the presentation as the maximum
filling area for words shorter than a given length.  That is,
$$\delta(n)=\max_{w\in R\cap B(n)}\delta(w),$$ where $B(n)$ denotes the ball of radius $n$ in $F$.
This depends on the choice of presentation, but changing the
presentation changes the asymptotics only minimally.  In particular,
if $\delta(n)=\Theta(n^k), k>1$ for one finite presentation, the
same is true for any other finite presentation of $G$.  If we view
$\delta(w)$ as a distance function on $R$ as above, then the Dehn
function of a group measures the distortion of $\delta(w)$ compared to
the metric induced on $R$ by inclusion in $F$

Instead of taking the maximum, however, we can average the filling
area over all words shorter than a given length that represent the
identity to get an averaged version, $\delta^{\text{avg}}$ of the Dehn function.  In order to
apply results on random walks, we will average over all
``lazy'' words of exactly a given length, that is, words of the form
$a_1\dots a_{n}$, where $a_j\in\{e_i^{\pm1},e\}$.

We will define $\delta^{\text{avg}}$ using random walks.
Let $p\in M(G)$ be the measure $$p(g)=\begin{cases}\frac{1}{2d+1} &
g=e \text{ or } g=e_i^{\pm1}\\ 0 & \text{otherwise.}\end{cases}$$ (For
our purposes, any finitely supported probability measure such that
$p(e)>0$, $p(x)=p(x^{-1})$, and the support of $p$ generates $G$ will suffice.) We use
$p$ to construct a random walk where at each step, the probability of
moving from $g$ to $gh$ is given by $p(h)$.  Then $p^{(n)}(x)$, the
$n$th convolution power of $p$, is the probability that an $n$-step random walk starting at $e$ ends at $x$.  We also define
$p^{(n)}(x,y)\equiv p^{(n)}(x^{-1}y)$, the probability of going from
$x$ to $y$ in $n$ steps and $p(x,y)\equiv p^{(1)}(x,y)$.

Define a measure $\rho_n$ on $G^n=G\times\dots\times G$ by
$$\rho_n(g_1,\dots,g_n)=\prod_{i=1}^{n-1} p(g_i).$$ Then
$\rho_n(g_1,\dots,g_n)$ represents the probability that an $n$-step
random walk starting at $e$ goes to $g_1$, then $g_1g_2$, and so on,
ending at $g_1\dots g_n$.  We consider this as a probability measure
on the set of lazy words of length $n$, where $(g_1,g_2,\dots,g_n)\in
G^n$ corresponds to the unreduced lazy word $g_1g_2\dots g_n$.

We can then define our averaged Dehn function by considering the
measure $\rho_{n}|_{\{(g_1,\dots,g_n)|g_1\dots g_n=e\}}$.  This
measure is nonzero, since $\rho_{n}(e,e,\dots,e)=p(e)^{n}>0$, so we
can normalize it to a probability measure $\overline{\rho}_{n}$ and
consider its support as the set of lazy words of length $n$ which are
the identity in $G$.  We then define the averaged Dehn function
$$\delta^{\text{avg}}(n)=E_{\overline{\rho}_{n}}(\delta(w))=\sum_w\delta(w)\overline{\rho}_{n}(w)$$
as the expected area necessary to fill a random closed curve of length $n$.
With the
choice of measure given above, the probability of any lazy word is the
same, and this average is the same as averaging over all lazy words
using the counting measure.

Like the Dehn function, the averaged Dehn function depends {\it a
priori} on the given presentation.  For nilpotent groups, we
will show the following upper bound, which is independent of the presentation:
\begin{thm}\label{upbound}
If $G$ is a finitely generated nilpotent group with Dehn function $\delta(n)=O(n^k)$, then
if $k>2$,
its averaged Dehn function for any presentation satisfies
$\delta^{\text{avg}}(n)=O(n^{k/2})$, and if $k=2$, $\delta^{\text{avg}}(n)=O(n \log n)$
\end{thm}
By a theorem of Gersten, Holt, and Riley\cite{GeHoRi}, all nilpotent
groups have Dehn functions bounded above by a polynomial, so as a
consequence, nilpotent groups with Dehn function growing at least
quadratically have averaged Dehn function strictly subasymptotic to
their Dehn function.  On the other hand, if a group has subquadratic
Dehn function, it is hyperbolic by a theorem of Gromov, one proof of which can be found in \cite{Bow}.  Since
nilpotent groups are amenable, a nilpotent hyperbolic group must be
finite or virtually $\Z$, so if $G$ is a nilpotent group which is not
finite or virtually $\Z$, its averaged Dehn function is strictly
subasymptotic to its Dehn function.

\section{Centralized Isoperimetry}
\label{centiso}
The lower bound we will give uses the centralized isoperimetric
function defined by Baumslag, Miller and Short \cite{BaMiSh}; we recall the definition.  Let $G=F/R=\{e_1,\dots,e_d|r_1,\dots,r_s\}$ as above.  
If $w$ is a word which is the
identity in $G$, we can write 
$$w\in\prod_{i=1}^k g_i^{-1} r_{a_i}^{\pm1} g_i[R,F]=\prod_{i=1}^k r_{a_i}^{\pm1}[R,F]=\prod_{i=1}^s r_i^{b_i}[R,F],$$
using the fact that $R[R,F]\subset Z(F/[R,F])$.
Define $\delta_G^{\text{cent}}(w)$ as the minimal $k$ for
which we can write such a decomposition (equivalently, the minimal $\sum_{i=1}^s|b_i|$) , and for $n\in\N$, define 
$$\delta_G^{\text{cent}}(n)=\max_{w\in R\cap B(n)}\delta_G^{\text{cent}}(w).$$  As before, we will drop the $G$ if the group is clear.  This depends {\it a priori} on the choice of presentation, but Baumslag, Miller, and Short \cite{BaMiSh} prove that, like the Dehn function, changing the
presentation changes the asymptotics only minimally.  In particular,
as for the Dehn function, if $\delta^{\text{cent}}(n)=\Theta(n^k),
k>1$ for one finite presentation, the same is true for any other
finite presentation of $G$.

Equivalently, give $F/[R,F]$ the generating set $\{e_1,\dots,e_d\}$
and give $R/[R,F]$ the generating set $\{r_1,\dots,r_s\}$.  If we
denote the distance functions induced by these generators by
$d_{F/[R,F]}$ and $d_{R/[R,F]}$, then

$$\delta^{\text{cent}}(w)=d_{R/[R,F]}(e,w)$$
and
$$\delta^{\text{cent}}(n)=\max_{w\in R/[R,F]\cap B_{F/[R,F]}(n)}d_{R/[R,F]}(e,w).$$
Thus, in the same way that the Dehn function measures the distortion
of the inclusion $R\subset F$ for the metric on $R$ induced by the
generating set $\{g^{-1}r_i g\}_{g\in G, 1\le i\le s}$, the
centralized isoperimetric function measures the distortion of the
inclusion $R/[R,F]\subset F/[R,F]$ for the metric on $R/[R,F]$ induced
by the generating set $\{r_i\}_{1\le i\le s}$.  Since
$R/[R,F]$ is finitely generated and abelian, $\delta^{\text{cent}}$ is
generally easier to calculate than $\delta$ and provides a lower bound
for it, since $$\delta^{\text{cent}}(w)\le \delta(w).$$

We can now state our lower bound on the averaged isoperimetric function:
\begin{thm}\label{lowbound}
  If $G$ is a finitely generated nilpotent group with centralized
  isoperimetric function $\delta^{\text{cent}}(n)=\Omega(n^k)$ for
  $k\ge 2$, then its averaged Dehn function for any presentation
  satisfies $\delta^{\text{avg}}(n)=\Omega(n^{k/2})$.
\end{thm}

Many nilpotent groups have $\delta^{\text{cent}}(n)$ and
$\delta(n)$ both polynomial of the same degree.  If in addition, this degree is $>2$, our upper and lower bounds are
sharp and independent of the presentation of the group.

In the remainder of this section, we will prove some results on
centralized isoperimetric functions of nilpotent groups which will be
useful in the proof of Theorem \ref{lowbound}.  If $G=F/R$ is a
finitely generated nilpotent group, $F/[R,F]$ is nilpotent and in fact
a central extension of $F/[R,F]$.  The asymptotics of
$\delta^{\text{cent}}$ are then relatively straightforward to
calculate.  Since $R/[R,F]\subset F/[R,F]$ is the inclusion of an
abelian group into a nilpotent group, we can apply a special case of a
theorem of Osin\cite{Osin}:

\begin{thm}
  Let $G$ be a f.g. nilpotent group, $H$ be an abelian subgroup of
  $G$, and $H^0$ the set of all elements of infinite order in $H$.  If
  $k$ is maximal such that $H^0\cap G^{(k)}\ne \{e\}$, then
  $$\max_{h\in H\cap B(n)}d_H(h,e)=\Theta(n^k).$$
  In particular, if $k$ is as above, $w\in H^0\cap G^{(k)}$ and $w\ne
  e$, then $d_H(w^{n^k},e)=\Theta(n).$
\end{thm}
where $G^{(\cdot)}$ is the lower central series of $G$,
\begin{align*}
G^{(1)}&=G,\\
G^{(n)}&=[G,G^{(n-1)}].
\end{align*}

We will call $k$ the {\em degree of distortion} of $H$ in $G$.  Then
\begin{align*}
\delta^{\text{cent}}(n)&=\max_{w\in R/[R,F]\cap B_{F/[R,F]}(n)}d_{R/[R,F]}(e,w).\\
    &=\Theta(n^k).
\end{align*}
where $k$ is the degree of distortion of $R/[R,F]$ in $F/[R,F]$.

In fact, we will show that any central extension of $G$ provides a
lower bound on its centralized isoperimetric function and that in fact
$\delta^{\text{cent}}(n)$ can be calculated by considering just
central extensions of $G$ by $\Z$.  These lower bounds are closely
related to the lower bounds for the centralized isometric function
found in Theorem 8 of \cite{BaMiSh} and in \cite{PitJLMS}.  We will
prove the following proposition, which we will use in proving a lower
bound on the averaged isoperimetric function:
\begin{prop} \label{prop:centext}
  If $G$ is a finitely generated nilpotent group, $A$ is a finitely
  generated abelian group,
  $$0\to A\to H \to G\to 1$$
  is a central extension of $G$, and $k$ is
  the degree of distortion of $A$ in $H$, then if $k\ge 2$(in
  particular, the extension must be nontrivial), then
  $\delta_G^{\text{cent}}(n)=\Omega(n^k)$.

  If $G$ is a finitely generated nilpotent group and
  $\delta_G^{\text{cent}}(n)=\Theta(n^k)$, then there is a central
  extension
  $$0\to \Z\to H \to G\to 1$$
  such that the degree of distortion of $\Z$ in $H$ is $k$.
\end{prop}
We prove the proposition by using the fact that a central extension by
$A$ can be described by a map $R/[R,F]\to A$.  Recall that any central extension of a nilpotent group is again nilpotent.
\begin{lemma}
  If $G=F/R=\{e_1,\dots,e_d|r_1,\dots,r_s\}$, $A$ is abelian and
  $$0\to A\stackrel{i}{\to} H \stackrel{p}{\to} G\to 1$$
  is a central extension of $G$ by $A$, then there are maps $\alpha:
  R/[R,F]\to A$ and $\beta: F/[R,F]\to H$ such that
  \begin{equation*}\label{extdiag}\xymatrix{
      0\ar[r]  & R/[R,F] \ar[r]\ar[d]^\alpha& F/[R,F] \ar[r]\ar[d]^\beta& G \ar[r]\ar[d]^\cong & 1\\
      0\ar[r]  & A \ar[r]& H \ar[r]& G \ar[r] & 1\\
    }\end{equation*}
  commutes and $\image\;\alpha \supset A \cap [H,H]$.

\end{lemma}
\begin{proof}[Proof of lemma]
  Choose $e'_1,\dots,e'_d\in H$ such that $p(e'_i)=e_i$.  This defines
  a map $f:F\to H$ which sends $e_i$ to $e'_i$.  We will find $\alpha$
  and $\beta$ from this map.  First, note that $f(R)\subset\ker p$, so
  $f(R)\subset A$.  Since $A$ is in the center of $H$,
  $f([R,F])\subset [A,H]=\{1\}$ and so we can define $\beta:F/[R,F]\to
  H$ as the quotient of $f$ by $[R,F]$ and $\alpha:R/[R,F]\to A$ as
  the restriction of $\beta$ to $R/[R,F]$.  These maps make the
  diagram commute.  Finally, if $\prod_{i=1}^n [h_{i,1},h_{i,2}]\in A
  \cap [H,H]$ for $h_{i,j}\in H$, choose $h'_{i,j}\in F$ such that
  $p(f(h'_{i,j}))=p(h_{i,j})$ (possible because $p\circ f$ is
  surjective).  Then
  $$p\left(f\left(\prod_{i=1}^n [h'_{i,1},h'_{i,2}]\right)\right)=p\left(\prod_{i=1}^n [h_{i,1},h_{i,2}]\right)=1$$
  and so $\prod_{i=1}^n [h'_{i,1},h'_{i,2}]\in R$.  Since
  $f(h'_{i,j})=h_{i,j}a_{i,j}$ for some $a_{i,j}\in A$,
  $$f\left(\prod_{i=1}^n[h'_{i,1},h'_{i,2}]\right)=\prod_{i=1}^n [h_{i,1}a_{i,1},h_{i,2}a_{i,2}]=\prod_{i=1}^n [h_{i,1},h_{i,2}]$$ 
  and so $\prod_{i=1}^n [h_{i,1},h_{i,2}]\in \image\;\alpha$ as desired.
\end{proof}

\begin{proof}[Proof of proposition]
  For the first part, let $\overline{A}=\image\;\alpha$.  Note that
  $\overline{A}\supset A \cap [H,H]$, so for all $j\ge 2$, we have
  $\overline{A}^0\cap H^{(j)}=A^0\cap H^{(j)}$.  Since the degree of
  distortion of $A$ in $H$ is $\ge 2$, it is the same as that of
  $\overline{A}$ in $H$ and so
  $$\max_{a\in \overline{A}\cap B_{H}(n)}d_{\overline{A}}(a,e)=\Theta(n^k).$$
  However, $\alpha$ cannot increase distances by more than a constant
  multiple, so $$d_{\overline{A}}(\alpha(w),e)\le C d_{R/[R,F]}(w,e),$$
  and
  \begin{align*}
    \delta^{\text{cent}}(n)&=\max_{w\in R/[R,F]\cap B_{F/[R,F]}(n)}d_{R/[R,F]}(e,w).\\
    &\ge \frac{1}{C}\max_{a\in \overline{A}\cap B_{H}(n)}d_{\overline{A}}(a,e)\\
    &=\Omega(n^k).
  \end{align*}
  When applied to extensions by $\Z$, this bound is a discrete
  analogue of the lower bound found by Pittet in \cite{PitJLMS}.

  For the second part, suppose that
  $\delta_G^{\text{cent}}(l)=\Omega(l^k)$.  We want to find a central
  extension by $\Z$ giving this bound.  Since the degree of distortion
  of $R/[R,F]$ in $F/[R,F]$ is at least $k$, there is an element $z$
  of $R/[R,F]$ of infinite order so that $z\in (F/[R,F])^{(k)}$.  Let
  $\alpha:R/[R,F]\to \Z$ be a map such that $\alpha(z)\ne 0$.  We can
  construct a central extension of $G$ by $\Z$ such that
\begin{equation*}\xymatrix{
0\ar[r]  & R/[R,F] \ar[r]\ar[d]^\alpha& F/[R,F] \ar[r]\ar[d]^\beta& G \ar[r]\ar[d]^\cong & 1\\
0\ar[r]  & \Z \ar[r]& H \ar[r]& G \ar[r] & 1\\
}\end{equation*}
commutes by letting $H=(F/[R,F])/\ker \alpha$.  Then $\alpha(z)\in \Z$ is of infinite order and $\alpha(z)\in H^{(k)}$, so
$$\max_{a\in \Z\cap B_{H}(n)}d_{\Z}(a,e)=\Omega(n^k).$$
Thus, to compute $\delta^{\text{cent}}$, it suffices to consider central extensions of $G$ by $\Z$.
\end{proof}

\section{Lower Bounds}\label{lowbndsec}
Here we extend the bounds in Section \ref{centiso} to the averaged
Dehn function.

We will be using theorems on the behavior of random walks on nilpotent groups, most notably a theorem of Hebisch and Saloff-Coste:
\begin{thm}[\cite{HeSal}]
\label{rndwalks}
Let $G$ be a finitely generated group with polynomial volume growth of
order $D$.  Let $p$ be a finitely supported probability measure such that
$p(e)>0$, $p(x)=p(x^{-1})$, and the support of $p$ generates $G$.  Then there exist three positive constants $C, C', C''$ such that, for all $x\in G$ and all integers $n$, we have
$$p^{(n)}(x)\le Cn^{-D/2}\exp\left(-d(e,x)^2/C'n\right)$$
$$p^{(n)}(x)\ge (Cn)^{-D/2}\exp\left(-C'd(e,x)^2/n\right) \quad\text{if $x\in B(e,n/C'')$}$$
where $d(e,x)$ denotes distance in the word metric corresponding to the support of $p$ and $B(e,r)$ is the ball of radius $r$ around $e$ in this metric.
\end{thm}

We will prove Theorem \ref{lowbound} by applying Theorem
\ref{rndwalks} to a suitable extension of $G$.  We first sketch an
outline of the proof, then fill in the calculational details.  By
Proposition \ref{prop:centext}, we know that if $G$ is a finitely
generated nilpotent group and $\delta_G^{\text{cent}}(n)=\Theta(n^k)$,
then there is a central extension
$$0\to \Z \to H \to G \to 1$$ 
such that
$$\max_{a\in \Z\cap B_{H}(n)}d_{\Z}(a,e)=\Omega(n^k).$$
If we choose $a_i'\in H$ that project to a set of generators $a_i\in G$,
then lazy words in the $a_i$ that are the identity in $G$ correspond
exactly to lazy words in the $a_i'$ representing elements in the image
of $\Z$.  Theorem \ref{rndwalks} suggests that a random lazy word of
length $n$ in a nilpotent group will represent an element of distance
on average $\sim\sqrt{n}$ from the identity.  Because of the order $k$
distortion of $\Z$ in $H$, this suggests that, on average, those lazy
words which are elements of $\Z$ will represent elements of distance
$\sim n^{k/2}$ from the identity in $\Z$.  Such a word $w$ will have
$\delta^{\text{cent}}(w)\approx n^{k/2}/c$.  Formalizing this argument
and carefully calculating the probabilities will give the lower bound
we want.

\begin{proof}[Proof of Theorem \ref{lowbound}]
  Construct $H$ and $a'_i$ as above.  One minor problem is that the
  $a'_i$ may not generate $H$, but unless $H=G\times \Z$(in which case
  $k=1$), the subgroup generated by the $a'_i$ is finite index in $H$
  and we can replace $H$ by this subgroup.

  Lazy words of length $n$ in $H$ will represent powers $z^l$ of $z$;
  we want to calculate the expected absolute value of $l$.
  Letting $p_H$ be the measure corresponding to the random walk on $H$
  with generating set $\{a'_1,\dots,a'_d\}$, this is
  $$\frac{\sum_{l\in\Z}|l|\cdot p_H^{(n)}(e,z^l)}{\sum_{l\in\Z}p_H^{(n)}(e,z^l)}.$$
  Note that since $G$ is nilpotent, it has polynomial volume growth,
  say of order $D$.  We thus apply Theorem \ref{rndwalks} to find
\begin{align*}
\frac{\displaystyle\sum_{l\in\Z}|l|\cdot p_H^{(n)}(e,z^l)}{\displaystyle\sum_{l\in\Z}{p_H^{(n)}(e,z^l)}}&
\ge \frac{\displaystyle\sum_{z^l\in B(e,n/C'')}|l|\cdot (Cn)^{-D/2}\exp\left(-C'd(e,z^l)^2/n\right)}{\displaystyle\sum_{l\in\Z}Cn^{-D/2}\exp\left(-d(e,z^l)^2/C'n\right)}\\
\end{align*}
For clarity, we will replace expressions not depending on $n$ with positive constants $c_i$.
\begin{align*}
\frac{\displaystyle\sum_{l\in\Z}|l|\cdot p_H^{(n)}(e,z^l)}{\displaystyle\sum_{l\in\Z}{p_H^{(n)}(e,z^l)}}
&\ge c_0\frac{\displaystyle\sum_{|l|<c_1n^k}|l|\cdot \exp\left(-c_2(|l|^{1/k})^2/n\right)}{\displaystyle \sum_{l\in\Z}\exp\left(-(|l|^{1/k})^2/c_2n\right)}\\
&\ge c_0\frac{\displaystyle\sum_{l=0}^{c_1n^k}l\cdot \exp\left(-c_2l^{2/k}/n\right)}{\displaystyle \sum_{l=0}^\infty \exp\left(-l^{2/k}/c_2n\right)}\\
&\ge c_0\frac{\displaystyle\int_0^{c_1n^k}l\cdot \exp\left(-c_2l^{2/k}/n\right)\;dl+O(n^{k/2})}{\displaystyle \int_0^\infty \exp\left(-l^{2/k}/c_2n\right)\;dl+1}
\end{align*}
We make the substitutions $x=(c_2/n)^{k/2}l$ and $y=(1/c_2n)^{k/2}l$.
\begin{align*}
\frac{\displaystyle\sum_{l\in\Z}|l|\cdot p_H^{(n)}(e,z^l)}{\displaystyle\sum_{l\in\Z}{p_H^{(n)}(e,z^l)}}
&\ge c_3\frac{\displaystyle n^k\int_0^{c_4n^{k/2}}x\cdot \exp\left(-x^{2/k}\right)\;dx+O(n^{k/2})}{\displaystyle n^{k/2} \int_0^\infty \exp\left(-y^{2/k}\right)\;dy+1}\\
&=\Omega(n^{k/2})
\end{align*}
\end{proof}

\section{Upper Bounds}\label{upbndsec}
To obtain upper bounds, we will construct discs filling in random
loops and bound their expected areas using Theorem \ref{rndwalks}.

Let $w=a_1a_2\dots a_n=e$ be a word with
$a_i\in\{e_1^{\pm1},\dots,e_d^{\pm1},e\}$.  We will think of $w$ as a
path in the Cayley graph and write $w(i)=a_1\dots a_i$.  Fix shortest
paths $\gamma_{x,y}$ between each pair of elements $x,y$ of $G$ so
that $\gamma_{x,y}$ is $\gamma_{y,x}$
traced backwards.  

\begin{figure}
  \centerline{\relabelbox \small \epsfbox{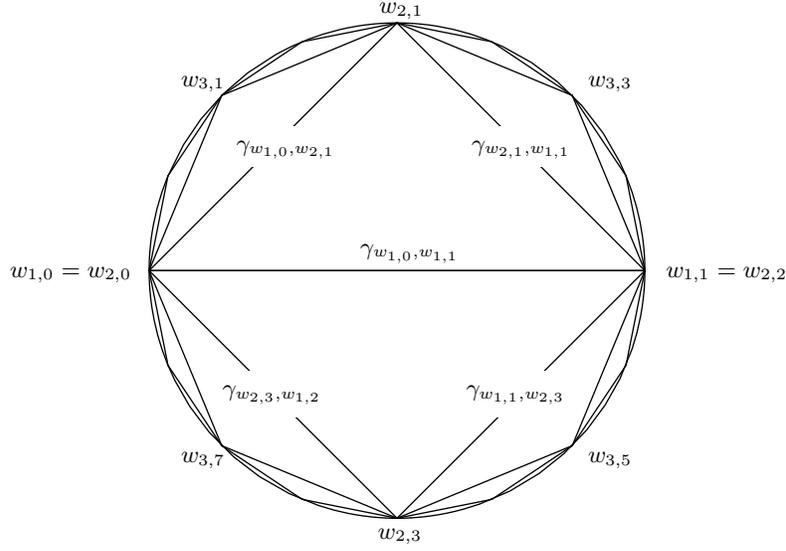} 
    \adjustrelabel <.12in,-.15in> {w1}{$w_{2,1}$}
    \relabel{w2}{$w_{1,1}=w_{2,2}$}
    \adjustrelabel <.17in,.17in> {w3}{$w_{2,3}$}
    \relabel{w4}{$w_{1,0}=w_{2,0}$}
    \relabel{g1}{$\gamma_{w_{1,0},w_{1,1}}$}
    \relabel{g2}{$\gamma_{w_{1,0},w_{2,1}}$}
    \relabel{g3}{$\gamma_{w_{2,1},w_{1,1}}$}
    \relabel{g4}{$\gamma_{w_{1,1},w_{2,3}}$}
    \relabel{g5}{$\gamma_{w_{2,3},w_{1,2}}$}
    \extralabel <-1in, .75in> {$w_{3,5}$}
    \extralabel <-3.13in, 2.72in> {$w_{3,1}$}
    \extralabel <-1in, 2.72in> {$w_{3,3}$}
    \extralabel <-3.13in, .75in> {$w_{3,7}$}
    \endrelabelbox }
  \caption{Filling a disc with triangles}
\label{discfill}
\end{figure}

Let
$$w_{i,j}=w\left(\left\lfloor \frac{j n}{2^{i}}\right\rfloor\right).$$
for $i\ge 1$, $0\le j\le 2^{i}$.  Note that $w_{i,j}=w_{i+1,2j}$ and
that for any $i$, $w_{i,1},\dots,w_{i,2^i}$ is an approximation of
$w$.  We will inductively build a sequence of fillings such that the
boundary of the $i$th filling is $w_{i+1,1},\dots,w_{i,2^{i+1}}$ by gluing triangles as in Figure \ref{discfill}.

We start with two discs filling
$$\gamma_{w_{1,0},w_{2,1}}\gamma_{w_{2,1},w_{1,1}}\gamma_{w_{1,1},w_{1,0}}$$
and
$$\gamma_{w_{1,1},w_{2,3}}\gamma_{w_{2,3},w_{1,2}}\gamma_{w_{1,2},w_{1,1}}$$
Gluing these discs together, we obtain a filling of 
$$\gamma_{w_{1,0},w_{2,1}}\gamma_{w_{2,1},w_{1,1}}\gamma_{w_{1,1},w_{2,3}}\gamma_{w_{2,3},w_{1,2}}=\gamma_{w_{2,0},w_{2,1}}\gamma_{w_{2,1},w_{2,2}}\gamma_{w_{2,2},w_{2,3}}\gamma_{w_{2,3},w_{2,4}}$$
After the $i$th step, we have a disc filling the geodesic $2^{i+1}$-gon with vertices
$$w_{i+1,1},\dots, w_{i+1, 2^{i+1}}.$$
To refine this to a filling of the next polygon, we add $2^{i+1}$ discs, filling the geodesic triangles with vertices 
$$w_{i+1,j}=w_{i+2,2j},w_{i+2,2j+1},w_{i+1,j+1}=w_{i+1,2j+2}.$$
Finally, after the $\lfloor \log_2 n\rfloor$th step, the boundary is
almost $w$.  We can apply a number of relators linear in $n$ to get
$w$ exactly.

It remains to estimate the total area of the triangles in the $\lfloor \log_2 n\rfloor$th filling.
\begin{proof}[Proof of Theorem \ref{upbound}]
  For a random closed path $w$ of length $n$ in $G$, we can construct
  a disc filling $w$ using the process above.  We must find a rigorous
  estimate of the area.  If $x, y, z\in G$, define $\Delta(x,y,z)$ to be
  the filling area of a geodesic triangle with vertices $x,y,z$, that is,
  $$\Delta(x,y,z)=\delta(\gamma_{x,y}\gamma_{y,z}\gamma_{z,x}).$$
  The process above gives a bound 
  $$
  \delta(w)\le cn +
  \sum_{i=1}^{\lfloor\log_2{n}\rfloor}\sum_{j=0}^{2^i-1}\Delta(w_{i,j},w_{i+1,2j+1},w_{i,j+1})$$
  Thus the expected area of a random word of length $n$ is at most
  \begin{align*}
    \delta^{\text{avg}}(n)&=E(\delta(w))\le E\left(cn +
  \sum_{i=1}^{\lfloor\log_2{n}\rfloor}\sum_{j=0}^{2^i-1}\Delta(w_{i,j},w_{i+1,2j+1},w_{i,j+1})\right)\\
    &= cn +
  \sum_{i=1}^{\lfloor\log_2{n}\rfloor}\sum_{j=0}^{2^i-1}E(\Delta(w_{i,j},w_{i+1,2j+1},w_{i,j+1})),
  \end{align*}
  where the expectations are taken with respect to
  $\overline{\rho}_{n}(\delta(w))$.  We can bound the averaged Dehn
  function by bounding the expected area of each triangle.  Consider 
  $E(\Delta(w(r),w(s),w(t)))$, the expected area of the triangle with
  vertices $w(r),w(s),$ and $w(t)$.  Since the Dehn function of $G$ is
  $O(n^k)$,
  $$\Delta(w(r),w(s),w(t))\le c_\delta (d(w(r),w(s))+d(w(s),w(t))+d(w(t),w(r)))^k.$$
  for some constant $c_\delta$, that is, the expected area is bounded
  by the $k$th moment of the perimeter.  This can be bounded with a
  classical inequality proved using the triangle inequality for $L^k$
  spaces; if we define $d_{s,t}(w)=d(w(s),w(t))$ as a function from
  the set of words of length $n$ representing the identity to $\R$ and
  give the set of such functions the $L^k$ norm $||\cdot||_k$ then
  $E(d_{s,t}(w)^k)=(||d_{s,t}||_k)^k$ and
  \begin{align}
    E(\Delta(w(r),w(s),w(t)))&\le E(c_\delta (d_{r,s}(w)+d_{s,t}(w)+d_{t,r}(w))^k)\notag\\
    &=(||d_{r,s}+d_{s,t}+d_{t,r}||_k)^k\notag\\
    &\le (||d_{r,s}||_k+||d_{s,t}||_k+||d_{t,r}||_k)^k\notag\\
    &\le 3^k(\max\{||d_{r,s}||_k,||d_{s,t}||_k,||d_{t,r}||_k\})^k\notag\\
    &\le 3^k\max\{E(d_{r,s}(w)^k),E(d_{s,t}(w)^k),E(d_{t,r}(w)^k)\}\label{triIneq}\tag{*}
  \end{align}
  Thus we can bound the expectation of
  $\Delta(w_{i,j},w_{i+1,2j+1},w_{i,j+1})$ by considering the $k$th
  moments of distances $d(w(s),w(t)).$

  We first claim that if $w$ is chosen from the unreduced words of
  length $n$ representing the identity in $G$ according to the
  probability distribution $\overline{\rho}_{n}$, then the
  distribution of $d(w(s),w(t))$ depends only on $|s-t|$ and $n$, that is,
  $$P(\{w|d(w(s),w(t))=x\})=P(\{w|d(w(0),w(t-s))=x\})$$
  This is true because the family of maps $r_i$ taking
  $$a_1\dots a_n$$ 
  to
  $$a_{i+1}a_{i+2}\dots a_na_1a_2\dots a_{i}$$ 
  preserves $\overline{\rho}_{n}$ and
  $$d(w(s),w(t))=d(e,r_s(w)(t-s)).$$ 
  Then
  \begin{align*}
    P(\{w|d(w(s),w(t))=x\}) &=\overline{\rho}_{n}(r_s^{-1}(\{w|d(e,w(t-s))=x\}))\\
    &=P(\{w|d(e,w(t-s))=x\})
  \end{align*}
  
  Thus to estimate the distribution of $d(w(s),w(t))$ it suffices to
  estimate the distribution of $d(e,w(t-s))$.  We will
  prove the
  following lemma:

  \begin{lemma}
    If $G$ is as above and $m\ge 1$, then there is a constant $c$ such
    that for any $n>0$, $t<n,$
    $$E\left(d(e,w(t))^m\right)<c t^{m/2},$$ 
    where the expectation is taken with respect to $\overline{\rho}_{n}$,
    i.e., over random closed paths of length $n$.
  \end{lemma}
  \begin{proof}
    Note that 
    $$d(e,w(t))=d(w(n),w(t)),$$
    so since the distribution of $d(w(s),w(t))$ depends only on
    $|s-t|$ and $n$, 
    $$E\left(d(e,w(t))^m\right)=E\left(d(e,w(n-t))^m\right),$$
    and we can assume that $t\le n/2$.  Note also that $d(e,w(t))\le
    t$, so $E\left(d(e,w(t))^m\right)\le t^m$ for all $t$.  It thus
    suffices to find bounds on $E\left(d(e,w(t))^m\right)$ for large
    $t$.
    
    The probability that a random closed path of length $n$ is at $x$ at time $t$ is:
    $$P(\{w|w(t)=x\})=\frac{p^{(t)}(e,x)p^{(n-t)}(x,e)}{\sum_{y\in G}p^{(t)}(e,y)p^{(n-t)}(y,e)}.$$
    Thus
    $$E\left(d(e,w(t))^m\right)=\frac{\sum_{x\in G}d(e,x)^m p^{(t)}(e,x)p^{(n-t)}(x,e)}{\sum_{y\in G}p^{(t)}(e,y)p^{(n-t)}(y,e)}$$
    Using Theorem \ref{rndwalks}, we can estimate these sums.   For
    clarity, we'll replace terms that don't depend on $n$ or $t$ with
    positive constants $c_i$.
    $$E\left(d(e,w(t))^m\right)\le c_0\frac{\sum_{x\in G}d(e,x)^m \exp\left(-d(e,x)^2/C't-d(e,x)^2/C'(n-t)\right)}{\sum_{y\in B(e,t/C'')}\exp\left(-C'd(e,y)^2/t-C'd(e,y)^2/(n-t)\right)}$$
    Since $0<t\le n/2$, $\frac{1}{t}\le\frac{1}{t}+\frac{1}{n-t}\le
    \frac{2}{t}$, so
    \begin{align*}
      E\left(d(e,w(t))^m\right)&\le c_0\frac{\sum_{x\in G}d(e,x)^m \exp\left(-\frac{c_1}{t} d(e,x)^2\right)}{\sum_{y\in B(e,t/C'')}\exp\left(-\frac{c_2}{t} d(e,y)^2\right)}\\
      &= c_0\frac{\sum_{r=0}^{\infty}\sum_{d(e,x)=r}r^m e^{-\frac{c_1}{t} r^2}}{\sum_{r=0}^{t/C''}\sum_{d(e,y)=r}e^{-\frac{c_2}{t} r^2}}\\
      &= c_0\frac{\sum_{r=0}^{\infty}\#S(e,r)r^m e^{-\frac{c_1}{t} r^2}}{\sum_{r=0}^{t/C''}\#S(e,r)e^{-\frac{c_2}{t} r^2}}
    \end{align*}
    where $\#S(e,r)$ is the number of $x\in G$ such that $d(e,x)=r$.
    We'd like to make the estimate $\#S(e,r)\approx r^{D-1}$, where
    $D$ is the order of polynomial growth of $G$, but this remains an
    open question.  Instead, we estimate the sums by a summation by
    parts argument that uses just the fact that $\#B(e,r)\approx r^D$.

    Recall that Abel's Formula states that if $\{a_i\}$ and
    $\{b_i\}$ are sequences and $B_n=\sum_{i=0}^n b_i$, then
    $$\sum_{i=n}^m a_i b_i=a_m B_m - a_{n-1}B_{n-1} - \sum_{i=n}^{m} B_i(a_{i+1}-a_{i})$$
    Let $\{d_i\}_{i\in \N}$ a sequence and $D_n=\sum_{i=1}^n d_i$ be
    its partial sums.  If $B_n \le D_n$ and $\{a_i\}$ is decreasing
    and positive, then
    \begin{align*}
      \sum_{i=n}^m a_i b_i&\le a_m D_m - a_{n-1}B_{n-1} - \sum_{i=n}^{m} D_i(a_{i+1}-a_{i})\\
      &= a_m D_m - a_{n-1}D_{n-1}+a_{n-1}(D_{n-1}-B_{n-1}) - \sum_{i=n}^{m} D_i(a_{i+1}-a_{i})\\
      &=\sum_{i=n}^m a_i d_i+a_{n-1}(D_{n-1}-B_{n-1})
    \end{align*}

    In particular, if $z n^D \le B_n\le Z n^D$ and $\{a_i\}$ is
    decreasing and positive, then
    $$\sum_{i=n}^m z' a_i i^{D-1}-a_{n-1}(Z-z)(n-1)^{D}\le\sum_{i=n}^m a_i b_i\le\sum_{i=n}^m Z' a_i i^{D-1}+a_{n-1}(Z-z)(n-1)^{D}$$
    for some $z',Z'$.

    To use this result, we need to replace $r^me^{-\frac{c_1}{t} r^2}$
    with a decreasing function of $r$.  The function has one extremum,
    a maximum of $\left(\frac{mt }{2c_1 e}\right)^{m/2}$ at
    $\sqrt{\frac{m t}{2c_1}}$.  Let $\beta_r=r^me^{-\frac{c_1}{t}
      r^2}$ and let $\beta'_r=e^{-\frac{c_2}{t} r^2}$.  Then we
    replace $\beta_r$ by $\left(\frac{m t}{2c_1 e}\right)^{m/2}$ when
    $r<\sqrt{\frac{m t}{2c_1}}$ to get:
    \begin{align*}
      E\left(d(e,w(t))^m\right)&\le c_0\frac{\sum_{r=0}^{\infty}\#S(e,r)\beta_r}{\sum_{r=0}^{t/C''}\#S(e,r)\beta'_r}\\
      &\le c_0\frac{\sum_{r=0}^{\sqrt{\frac{m t}{2c_1}}}\#S(e,r)\left(\frac{m t}{2c_1
            e}\right)^{m/2} +\sum_{r=\sqrt{\frac{m t}{2c_1}}}^{\infty}\#S(e,r)\beta_r}{\sum_{r=0}^{t/C''}\#S(e,r)\beta'_r}\\
&\le
      c_0\frac{\#B\left(e,\sqrt{\frac{m t}{2c_1}}\right)\left(\frac{m t}{2c_1
            e}\right)^{m/2}+\sum_{r=\sqrt{\frac{m t}{2c_1}}}^{\infty}\#S(e,r)\beta_r}{\sum_{r=0}^{t/C''}\#S(e,r)\beta'_r}
    \end{align*}

    Now $\beta_i$ and $\beta'_i$ are both decreasing in the intervals
    of summation.  Nilpotent groups have polynomial volume growth, so
    we can bound $\#B(e,r)$ by
    $$c_3^{-1} r^D\le\#B(e,r)\le c_3r^D\quad\text{for all $r\ge 1$.}$$
    Since $\#B(e,r)=\sum_{i=0}^r \#S(e,i)$, we can use the argument
    above to replace the $\#S(e,r)$'s by terms growing like $r^{D-1}$:
    \begin{align*}
      E\left(d(e,w(t))^m\right)&\le c_0\frac{c_3 \left(\sqrt{\frac{m t}{2c_1}}\right)^D \left(\frac{m t}{2c_1 e}\right)^{m/2}+\sum_{r=\sqrt{\frac{m t}{2c_1}}}^\infty c_4r^{D-1}r^me^{- \frac{c_1}{t} r^2}} {\sum_{r=0}^{t/C''}c_6 r^{D-1}e^{-\frac{c_2}{t} r^2}} \\
      &\le c_0\frac{c_5 t^{(D+m)/2}+\int_{\sqrt{\frac{m t}{2c_1}}}^\infty c_4r^{m+D-1}e^{-\frac{c_1}{t} r^2}\;dr+O(t^{(D+m-1)/2})} {\int_{0}^{t/C''}c_6 r^{D-1}e^{-\frac{c_2}{t} r^2}\;dr+O(t^{(D+1)/2})} 
    \end{align*}
    where in replacing the sums with integrals, we again use that
    $r^\lambda e^{-\frac{c}{t} r^2}$ has a maximum of $\left(\frac{\lambda t }{2c e}\right)^{\lambda/2}$

    If we let $\mu=D+m$ and substitute $x=\sqrt{\frac{1}{t}} r$, we get
    $$E\left(d(e,w(t))^m\right)\le c_0\frac{c_5 t^{\mu/2}+t^{\mu/2}\int_{\sqrt{\frac{m}{2c_1}}}^\infty c_4x^{\mu-1}e^{-c_1 x^2}dx+O(t^{(\mu-1)/2})}{t^{D/2}\int_{0}^{\sqrt{t}/C''}c_6 x^{D-1}e^{-c_2 x^2}dx+O(t^{(D-1)/2})} $$
    and thus
    $$E\left(d(e,w(t))^m\right)< c t^{m/2}$$
    as desired.
  \end{proof}

  One can use the same techniques to show that 
  $$E\left(d(e,w(t))^m\right)=\Theta(t^{m/2}).$$

  We will use this to bound the area of a triangle in the
  construction.  The triangles added in the $i$th step of the
  construction connect points $w(\lfloor j n/2^{i}\rfloor),
  w(\lfloor (2j+1) n/2^{i+1}\rfloor)$, and $w(\lfloor (j+1)
  n/2^{i}\rfloor)$.  By the lemma above, the $k$th moment of the
  expected distance between any two of these points is at most
  $c(n/2^{i})^{k/2}$, and thus by \eqref{triIneq},
  $$E(\Delta(w_{i,j},w_{i+1,2j+1},w_{i,j+1}))\le 3^kcc_\delta(n/2^{i})^{k/2}=Cn^{k/2}2^{-ki/2}.$$

  Finally, we find that
  \begin{align*}
    \delta^{\text{avg}}(n)    &= cn +
  \sum_{i=1}^{\lfloor\log_2{n}\rfloor}\sum_{j=0}^{2^i-1}E(\Delta(w_{i,j},w_{i+1,2j+1},w_{i,j+1}))\\
 &\le cn +
  \sum_{i=1}^{\lfloor\log_2{n}\rfloor}\sum_{j=0}^{2^i-1}Cn^{k/2}2^{-ki/2} \\
 &\le cn +
  \sum_{i=1}^{\lfloor\log_2{n}\rfloor}Cn^{k/2}2^{i(1-k/2)}.
  \end{align*}
  If $k>2$, this is a geometric series and
  $$\delta^{\text{avg}}(n)\le cn +
  \frac{C}{1-2^{1-k/2}} n^{k/2}=O(n^{k/2}).$$
  If $k=2$, then 
  $$\delta^{\text{avg}}(n)\le cn +
  \sum_{i=1}^{\lfloor\log_2{n}\rfloor}Cn=O(n \log n)$$
  as desired.
\end{proof}
For many nilpotent groups, $\delta(l)$ and $\delta^{\text{cent}}(l)$ have
polynomial growth of the same order; some examples of these are
abelian groups, free nilpotent groups\cite{BaMiSh}, and the
Heisenberg groups\cite{Allcock,BaMiSh}, though there are many more
examples.  For example, if $G$ is a nilpotent group of nilpotency
class $c$ such that $G^{(c)}$ contains elements of infinite order,
then $G/G^{(c)}$ has a Dehn function and a central isoperimetric function both with polynomial growth of order $c$.  The upper bound on the Dehn
function follows from \cite{GeHoRi}, and the lower bound on the
central isoperimetric function follows from Theorem 8 in \cite{BaMiSh}.
When this growth is faster than quadratic, as in non-abelian free
nilpotent groups, Theorems
\ref{upbound} and \ref{lowbound} give sharp estimates, independent of the generating set, of the growth of $\delta^{\text{avg}}$.

\section{Conclusion}
One natural question is how well these results extend to other groups.
The general idea that the loop generated by a random walk will have a
length scale much smaller than the number of steps taken seems likely
to hold in other groups, but the proofs here rely on good upper and
lower bounds for the off-diagonal transition probabilities of a random
walk, which may not be obtainable in other classes of groups.

One can consider the behavior of random closed paths in general.
For any $n$, we can construct a time-dependent random walk
$\hat{p}$ describing the behavior of random closed paths of length
$n$.  Assume that $p$ is symmetric, so that $p(x,y)=p(y,x)$.  The probability that a random
path of length $n$ from $x$ to $y$ is at $z$ after time $t$ is
$$\frac{p^{(t)}(x,z)p^{(n-t)}(z,y)}{\sum_w p^{(t)}(x,w)p^{(n-t)}(w,y)}.$$ 
Thus, if a random loop starting at $e$ is at $x$ after $t$ steps, it
must return to $e$ after $n-t$ more steps and we can write
$$\hat{p}(x,y;t)=\frac{p(x,y)p^{(n-(t+1))}(y,e)}{\sum_w p(x,w)p^{(n-(t+1))}(w,e)}$$
as the probability that its next step will take it to $y$.

We can often write 
$$\lim_{n\to\infty}\frac{p^{(n)}(e,x)}{p^{(n)}(e,e)}=f(x);$$
such a theorem is called a {\it ratio limit theorem}.  
In this case, 
\begin{align*}
\lim_{n\to\infty}\hat{p}(x,y;t)&=\lim_{n\to\infty}\frac{p(x,y)p^{(n-(t+1))}(y,e)}{\sum_w p(x,w)p^{(n-(t+1))}(w,e)}\\
      &=\frac{p(x,y)f(y)}{\sum_w p(x,w)f(w)}
\end{align*}
In amenable groups \cite{amenrat}, for instance, $f(x)=1$ for all $x$, and thus
$$\lim_{n\to\infty}\hat{p}(x,y;t)=p(x,y).$$
That is, when $n$ is large, a random closed path will look like the
standard random walk on small timescales.  The solvable
Baumslag-Solitar groups $BS(1,n)$ are examples of amenable groups with
exponential Dehn function for which the random walk goes to
infinity at a sublinear rate\cite{bsrwalks}, so it seems likely that
its averaged Dehn function is subexponential.

Sublinear growth of distances in random closed paths may be a fairly
general phenomenon.  In a free group, for example, the map $w\mapsto
d(w(i),e)$ taking random closed paths on the free group to random
closed paths on $\N$ is measure-preserving, so, as in the nilpotent
case, we find that $E(d(s),d(t))=O(\sqrt{|s-t|})$.

\def\cprime{$'$}

\end{document}